\documentclass[11pt,reqno]{amsart}
\usepackage{a4wide}
\usepackage{euscript,amsmath,amssymb,amsbsy,mathabx}
\usepackage[arrow,matrix,curve]{xy}
\usepackage{array,longtable,graphicx,enumitem,url}

\usepackage[colorlinks=true,linkcolor=black,urlcolor=black,citecolor=black]{hyperref}

\numberwithin{equation}{section}\numberwithin{figure}{section}

\newcounter{msct}[section]\renewcommand{\themsct}{\thesection.\arabic{msct}}
\newenvironment{m-theorem}{\vskip5pt\refstepcounter{msct}\trivlist \itemindent 0pt%
\item[\hskip\labelsep\bf Theorem~\themsct]\it\ignorespaces}{\endtrivlist\vskip3pt}
\newenvironment{m-proposition}{\vskip5pt\refstepcounter{msct}\trivlist \itemindent0pt%
\item[\hskip\labelsep\bf Proposition~\themsct]\it\ignorespaces}{\endtrivlist\vskip3pt}
\newenvironment{m-corollary}{\vskip5pt\refstepcounter{msct}\trivlist \itemindent 0pt%
\item[\hskip\labelsep\bf Corollary~\themsct]\it\ignorespaces}{\endtrivlist\vskip3pt}
\newenvironment{m-lemma}{\vskip5pt\refstepcounter{msct}\trivlist \itemindent 0pt%
\item[\hskip\labelsep\bf Lemma~\themsct]\it\ignorespaces}{\endtrivlist\vskip3pt}
\newenvironment{m-definition}{\vskip5pt\refstepcounter{msct}\trivlist \itemindent0pt%
\item[\hskip\labelsep\bf Definition~\themsct]\ignorespaces}{\endtrivlist\vskip5pt}
\newenvironment{m-notation}{\vskip5pt\refstepcounter{msct}\trivlist \itemindent0pt%
\item[\hskip\labelsep\bf Notation~\themsct]\ignorespaces}{\endtrivlist\vskip5pt}
\newenvironment{m-example}{\vskip5pt\refstepcounter{msct}\trivlist \itemindent0pt%
\item[\hskip\labelsep\bf Example~\themsct]\ignorespaces}{\endtrivlist\vskip5pt}
\newenvironment{m-remark}{\vskip5pt\refstepcounter{msct}\trivlist \itemindent0pt%
\item[\hskip\labelsep\bf Remark~\themsct]\ignorespaces}{\endtrivlist\vskip5pt}
\newenvironment{m-question}{\vskip3pt\refstepcounter{msct}\trivlist \itemindent0pt%
\item[\hskip\labelsep\bf Question.]\ignorespaces}{\endtrivlist\vskip5pt}
\newenvironment{thm-nono}[1]{\vskip5pt\trivlist \itemindent 0pt %
\item[\hskip\labelsep\bf Theorem~{\rm\mbox{#1}}]\it\ignorespaces}{\endtrivlist\vskip5pt}
\newenvironment{lm-nono}[1]{\vskip5pt\trivlist \itemindent0pt%
\item[\hskip\labelsep\bf Lemma~{\rm\mbox{#1}}]\it\ignorespaces}{\endtrivlist\vskip5pt}
\newenvironment{conj-nono}[1]{\vskip5pt\trivlist \itemindent0pt%
\item[\hskip\labelsep\bf Conjecture~{\rm\mbox{#1}}]\it\ignorespaces}{\endtrivlist\vskip5pt}
\newenvironment{m-thank}{\vskip5pt\trivlist \itemindent0pt%
\item[\hskip\labelsep\it Acknowledgments]\ignorespaces}{\endtrivlist\vskip5pt}
\newenvironment{m-proof}{\vskip2pt\trivlist \itemindent0pt%
\item[\hskip\labelsep\it Proof.]\ignorespaces}{\hfill$\Box$\endtrivlist\vskip5pt}%
\newenvironment{m-asmp}{\vskip5pt\trivlist \itemindent0pt%
\item[\hskip\labelsep\bf Assumption.]\ignorespaces}{\hfill\endtrivlist\vskip5pt}%

\newcommand{\bibauth}[2]{\textrm{{#1}~{#2}}}
\newcommand{\bibtitl}[1]{\textit{#1},}
\newcommand{\bibjnyp}[4]{\textrm{#1} \textbf{#2} (#3), #4.}
\newcommand{\bibinbook}[4]{In: \textrm{#1}\textrm{, #2}\textrm{, #3}\textrm{, #4}.}
\newcommand{\bibbook}[4]{\textit{#1}. {#2} {#3}, {#4}.}

\let\lar\longrightarrow

\let\mt\mapsto

\font\tenmsa=msam10 %
\newcommand\hdashpiece{%
{\vrule height2.75pt depth-2.35pt width2.3pt \kern1.7pt}}%
\newcommand\hdashpieces{%
{\hdashpiece\hdashpiece\hdashpiece\hdashpiece}}%
\let\dashto\dashrightarrow
\newcommand\dashar{\mathrel{%
\hdashpieces\kern-0.4pt\hbox{\tenmsa K}}}%

\let\euf\EuScript 
\let\cal\mathcal
\let\mbb\mathbb

\let\bsymb\boldsymbol 
\DeclareFontFamily{OT1}{rsfs}{}
\DeclareFontShape{OT1}{rsfs}{n}{it}{<->rsfs10}{}
\DeclareMathAlphabet{\crl}{OT1}{rsfs}{n}{it}

\newcommand\uset[2]{{\disp\mathop{\mbox{$#2$}}_{#1}}}
\newcommand\oset[2]{{\disp\mathop{\mbox{$#2$}}^{#1}}}
\newcommand\ouset[3]{{\oset{#2}{\uset{#1}{#3}}}}
\let\ovl\overline

\let\unbar\underbar

\let\tld\tilde
\let\wtld\widetilde

\let\nit\noindent

\let\disp\displaystyle
\let\srel\stackrel
\let\vphi\varphi

\let\veps\varepsilon

\newcommand\lran[1]{{\langle #1\rangle}}

\newcommand\ort{\mathrel{{\vrule width4.0pt height0.4pt depth0pt
                \vrule width0.4pt height6.0pt depth0pt\,}}}

\newcommand\Aut{\operatorname{\textrm{Aut}\kern1pt}}
\newcommand\cAut{\operatorname{\mathcal{A}\kern-1pt\textit{ut}\kern1pt}}
\newcommand\End{\operatorname{\rm{End}\kern1pt}}
\newcommand\cEnd{\operatorname{\mathcal{E}\kern-1pt\textit{nd}\kern1pt}}

\newcommand\Hom{\mathop{\rm Hom}\nolimits}
\newcommand\cHom{\operatorname{\mathcal{H}\kern-1pt\textit{om}\kern1pt}}
\newcommand\Img{{\rm Im}}
\newcommand\Ker{{\rm Ker}}

\newcommand\Pic{\mathop{\rm Pic}\nolimits}
\newcommand\Proj{\mathop{\rm Proj}\nolimits}

\newcommand\Spec{\mathop{\rm Spec}\nolimits}

\newcommand\Sym{\mathop{\rm Sym}\nolimits}
\newcommand\invq{{\slash\kern-2.5pt\slash}}

\newcommand\Grs{{\rm Gr}}
\newcommand\spGrs{{\rm sp{\hdashpiece\!}Gr}}
\newcommand\oGrs{{\rm o{\hdashpiece\!}Gr}}

\let\lda\lambda

\let\si\sigma

\let\sm\setminus

\newcommand\bbk{\mbox{\rm I\kern-1.5pt k}}
\newcommand\sbbk{\hbox{\scriptsize I{\kern-.8pt}k}}

\newcommand\bbQ{{\mbb Q}}

\newcommand\bone{{1\kern-0.57ex\rm l}}

\newcommand\eA{{\euf A}}

\newcommand\eE{{\euf E}}

\newcommand\eF{{\euf F}}
\newcommand\eG{{\euf G}}
\newcommand\eI{{\euf I}}
\newcommand\eJ{{\euf J}}

\newcommand\eL{{\euf L}}

\newcommand{\eM}{{\euf M}}

\newcommand\eN{{\euf N}}
\newcommand\cO{{\cal O}}
\newcommand\eO{{\euf O}}

\newcommand\eT{{\euf T}}

\newcommand\codim{{\rm codim}}

\newcommand\Bl{{\rm Bl}}

\newcommand\SO{{\rm SO}}

\newcommand\Sp{{\rm Sp}}

\renewcommand\det{{\rm det}}

\let\ges\geqslant
\let\les\leqslant

\newcommand\mn{{\rm min}}
\newcommand\res{\mathop{\rm res}\nolimits}
\let\surj\twoheadrightarrow

\newcommand\pos{{>0}}
\newcommand\apos{{\,\gtrsim0}}

\newcommand\diag{\mathop{\rm diag}\nolimits}

\newcommand{\BB}{{\rm BB}}

\newcommand{\cd}{\mathop{\rm cd}\nolimits}

\newcommand{\cst}[1]{\mathop{\rm ct}^{#1}\nolimits}

\newcommand{\kk}{{\Bbbk}}
\newcommand{\lcit}{{\textit{loc.\,cit.}}}
\newcommand{\ocit}{{\textit{op.\,cit.}}}
\newcommand{\red}{{\rm red}}
\newcommand{\reg}{\mathop{\rm reg}\nolimits}

\newcommand{\rst}{{\upharpoonright}}
\newcommand{\bz}{{\textbf{z}}}
\newcommand{\lbkt}{{[\kern-1.75pt[}}
\newcommand{\rbkt}{{]\kern-1.75pt]}}

\begin{document}

\title[Partially ample subvarieties]{Partially ample subvarieties of projective varieties}
\author{Mihai Halic}
\email{mihai.halic@gmail.com}
\address{CRM, UMI 3457, Montr\'eal~H3T~1J4, Canada}


\begin{abstract}
The goal of this article is to define partially ample subvarieties of projective varieties, generalizing Ottem's work on ample subvarieties, and also to show their ubiquity. As an application, we obtain a connectedness result for pre-images of subvarieties by morphisms, reminiscent to a problem posed by Fulton-Hansen in the late 1970s. Similar criteria are not available in the literature.
\end{abstract}

\maketitle


\section*{Introduction}

Hartshorne~\cite{hart-cdav,hart-as} pioneered the systematic study of the cohomological properties of pairs consisting of a projective scheme which is regular in a neighbourhood of a local complete intersection subscheme with ample normal bundle. Ample subvarieties of projective varieties were defined by Ottem~\cite{ottm}, based on Totaro's work~\cite{tot} on cohomological ampleness. Inspired by \ocit, but in a totally different framework (cf.~\cite[Introduction]{hlc-subvar}), the author considered the weaker $q$-ampleness property. Moreover, recently there is an increased interest in studying and understanding various positivity properties of higher codimensional subvarieties and cycles, see \textit{e.g.}~\cite{ful+leh} and the references therein. 

Our goal is to introduce partially ample subschemes. Besides general properties and examples, an important feature are their connectedness properties. The main result is a contribution to an issue raised by Fulton-Hansen~\cite{fult+hans}, known to be false in its original form. We prove that appropriate partial ampleness of a subvariety yields the connectedness of its pre-image under a morphism; in fact, it implies the G3 property, in Hironaka-Matsumura's terminology. 

The article consists of three sections. In the first one, we introduce the relevant definitions and \emph{compress} those properties which carry over from~\cite{ottm}; for details, the reader should consult the original reference. Next we present a class of situations where Fulton-Hansen's question does admit a positive answer.

\begin{thm-nono}{(cf.~\ref{thm:f-1},~\ref{thm:q})}
Let $V,X$ be irreducible projective varieties, $V\srel{f}{\to} X$ a morphism. Let $Y\subset X$ be a closed subscheme. 
\begin{enumerate}[leftmargin=5ex]
\item[\rm(i)] 
If $f$ is surjective and $\cd(X\sm Y)\les\dim X-2$, then $f^{-1}(Y)$ is connected. 
\item[\rm(ii)] 
Suppose $Y$ is $\big(\dim f(V)+\dim(Y)-\dim(X)-1\big)$-ample. 
\item[] 
Then $f^{-1}(Y)$ is connected and $\pi_1^{alg}(f^{-1}(Y))\to\pi_1^{alg}(V)$ is surjective.
\end{enumerate}
\end{thm-nono}
When $f$ is an embedding, the theorem yields a connectedness criterion for intersections, reminiscent to an problem posed by Hartshorne (cf.~\cite{hart-as,petr}) which is still wide open. An example from \ocit\, shows that the first part of our result is optimal (cf. Remark~\ref{rmk:optim}). To our knowledge, such numerical conditions are not available in this generality. Existing results (cf. Fulton-Hansen~\cite{fult+hans,hans}, Faltings~\cite{falt-homog}, Debarre~\cite{debr2}) hold for subvarieties of various homogeneous spaces. 

Also, \emph{our result goes beyond the applicability of Ottem's work}. On one hand, for a subvariety, the requirement to be partially ample is obviously less restrictive than to be ample. On the other hand---at a deeper level---even in the basic case of two subvarieties of an ambient space, one of them being ample, \cite{ottm} yields only the non-emptiness of their intersection, but gives \emph{no information about the connectedness} without further smoothness assumptions. Since the image of a morphism can have arbitrarily bad singularities, it is not possible to conclude anything regarding the connectedness of pre-images. Precisely for this reason, the first part of the Theorem is essential: it holds in full generality.

In the last section, we show the ubiquity of partially ample---\emph{mostly not ample}---subvarieties by analysing several classes of examples: vanishing loci of sections in vector bundles (cf. {\S}\ref{ssct:glob-gen}); Bialynicki-Birula decompositions (cf. {\S}\ref{ssct:fixed}); rational homogeneous spaces (cf. {\S}\ref{ssct:subvar-homog}).


\section{{\textit{q}-ample and \textit{p}-positive subvarieties}}

\begin{m-notation}\label{not:XYN} 
Let $X$ be a projective scheme defined over an algebraically closed field $\kk$ of characteristic zero. Let $Y$ be a closed subscheme; we denote the maximal dimension of its components by $\dim Y$, and assume that they are all at least $1$-dimensional. Let $\eI_Y\subset\eO_X$ be the sheaf of ideals defining $Y$; for $m\ges0$, $Y_m$ is the subscheme defined by $\eI_{Y}^{m+1}$. The formal completion of $X$ along $Y$ is $\hat X_Y:=\disp\varinjlim Y_m$; for any coherent sheaf $\eG$ on $X$, it holds 
\begin{equation}\label{eq:XY}
H^t(\hat X_Y,\eG)=\varprojlim H^t(Y_m,\eG). 
\end{equation}
Let $\tld X:=\Bl_Y(X)\srel{\si}{\to}X$ be the blow-up of $\eI_Y$ and $E_Y\subset\tld X$ the exceptional divisor. If $X$ is Cohen-Macaulay and $Y$ is locally complete intersection---\emph{lci} for short---, its normal sheaf $\eN_{Y/X}:=(\eI_Y/\eI_Y^2)^\vee$ is locally free.  A variety is a reduced and irreducible scheme. The symbol $\cst{A}$ stands for a constant depending on the quantity $A$. Further necessary notions are recalled in the appendix.
\end{m-notation}


\subsection{Basic properties}\label{sct:1-property}

Let the situation be as above.

\begin{m-definition}\label{def:q1}
We denote 
$$
\delta:=\codim_X(Y)=\min\{\codim_XY'\mid Y'\;\text{irreducible component of}\;Y\}.
$$ 
\begin{enumerate}[leftmargin=5ex]
\item[\rm(i)] {\rm(cf. \cite[Definition 3.1]{ottm})}  
We say that $Y$ is \emph{$q$-ample} if $\eO_{\tld X}(E_Y)$ is $(q+\delta-1)$-ample. That is, for any coherent sheaf $\tld\eF$ on $\tld X$ it holds: 
\begin{equation}\label{eq:q11}
H^{t}(\tld X,\tld\eF\otimes\eO_{\tld X}(mE_Y))=0,\;\forall\,t\ges q+\delta,\;\forall\,m\ges\cst{\tld\eF}.
\end{equation}
\item[\rm(ii)] 
We say that $Y$ is (has the property) $p^\pos$---that is, \emph{$p$-positive}---if it holds: 
\begin{equation}\label{eq:q12}
H^t(X,\eF\otimes\eI_Y^m)=0,\;\forall\,t\les p,\;\forall\,m\ges\cst{\eF},
\end{equation}
for all locally free sheaves $\eF$ on $X$.
\end{enumerate}
\end{m-definition}
Partial ampleness behaves well under restrictions, while positivity yields connectedness results. The notions are dual under certain regularity assumptions.
\begin{m-remark}
\nit{\rm(i)} For $q=0$, one recovers the ample subschemes~\cite{ottm}. 

\nit{\rm(ii)} For a closed point $y\in Y$, it holds 
\begin{equation}\label{eq:les}
\delta_{Y,y}-1\les\dim\si^{-1}(y)\les q+\delta-1,
\end{equation}
where $\delta_{Y,y}$ is the codimension of the irreducible component of $Y$ containing $y$ (cf.~\cite[Proposition 3.4]{ottm}). See~\ref{lm:conn} for equidimensionality criteria of partially ample subvarieties.
\end{m-remark}

\begin{m-proposition}\label{prop:N+cd} 
\begin{enumerate}[leftmargin=5ex]
\item[\rm(i)] 
A subscheme is $q$-ample if and only if so is its integral closure.
\item[\rm(ii)] 
The $q$-ampleness property of lci subschemes is open relative to projective, flat morphisms.
\item[\rm(iii)] 
$Y$ is $q$-ample if and only if it holds:
\item[]  
{\rm(a)} $\eO_{E_Y}(E_Y)$ is $(\delta+q-1)$-ample; 
{\rm(b)} $\cd(X\sm Y){\les}\,\delta+q-1.$ 
(Note:$\;\cd(X\sm Y){\ges}\,\delta-1$.)
\end{enumerate}
\end{m-proposition}

\begin{m-proof}
(i)+(ii) See \cite[Proposition~6.8, Theorem~6.1]{ottm}. 

\nit(iii) 
Verify~\eqref{eq:q11} for $\tld\eF_{E_Y}$---$\tld\eF={\tld\eA}^{-k}$, $k\ges 1$, and $\tld\eA\in\Pic(\tld X)$ ample---by using the sequence 
$0\to\tld\eF((m-1)E_Y)\to\tld\eF(mE_Y)\to\tld\eF_{E_Y}(mE_Y)\to0;$ 
proceed as in~\cite[Theorem 5.4]{ottm}.
\end{m-proof}

\begin{m-proposition}\label{prop:p>0}
The following statements are equivalent:
\begin{enumerate}[leftmargin=5ex]
\item[\rm(i)] $Y\subset X$ is $p^\pos$;\qquad {\rm(ii)} $E_Y\subset\tld X$ is $p^\pos$; 
\item[\rm(iii)] 
The condition~\eqref{eq:q11} holds for $\eF=\eA^{k}$, $k\ges1$, where $\eA\in\Pic(X)$ is ample. 
\item[\rm(iv)] 
For all locally free sheaves $\eF$ on $X$, the properties below are satisfied:\\[1ex] 
$\begin{array}{rl}
{\rm(a)}&
  H^t\big(X,\eI_Y^m/\eI_Y^{m+1}\otimes\eF\big)=0,\;\forall\,m\ges\cst{\eF},\;0\les t\les p-1,
\\[0ex]
{\rm(b)}&
  \res^X_{Y}:H^t(X,\eF)\to H^t(\hat X_Y,\eF) 
  \text{ is }\;\biggl\{
    \begin{array}{l}
    \text{an isomorphism, for }t\les p-1,\\[1ex] 
    \text{injective, for }t=p.
    \end{array}
  \Big.
\end{array}$
\end{enumerate}
\end{m-proposition}

\begin{m-proof}  
(i)$\Leftrightarrow$(iii) Observe that $\eF$ fits into $0{\to}\eF{\to}\eA^k\otimes\kk^N{\to}\eG{\to}0$, with $k,N>0$, and $\eG$ locally free. Then $H^t(\eF\otimes\eI_Y^m)\cong H^{t-1}(\eG\otimes\eI_Y^m)$---so $t$ decreases---and we repeat the process. 

\nit(i)$\Leftrightarrow$(ii) Let $\eA\in\Pic(X)$ be ample such that $\tld\eA:=\si^*\eA(-E_Y)$ is ample on $\tld X$. Now note that 
$\,H^t(\tld X,\tld\eA^{k}(-mE_Y)){=}
H^t(\tld X,\si^*\eA^{k}(-(k+m)E_Y)){=}
H^t(X,\eA^k\otimes\eI_Y^{k+m})$, for $m{\gg}0$.

\nit(i)$\Leftrightarrow$(iv)  For (a), we use $0{\to}\eI_Y^{m+1}{\to}\eI_Y^m{\to}\eI_Y^m/\eI_Y^{m+1}{\to}0$; for (b), twist by $\eF$ the following sequence and use~\eqref{eq:XY}: $0{\to}\eI_Y^{m+1}{\to}\eO_X{\to}\eO_{Y_m}{\to}0$. 

Conversely, for $t\les p-1$ and $m\gg\cst{\eF}$, using  $0{\to}\eI_Y^m/\eI_Y^{m+1}{\to}\eO_{Y_m}{\to}\eO_{Y_{m-1}}{\to}0$ we deduce that $H^t(Y_m,\eF){\to} H^t(Y_{m-1},\eF)$ are injective and eventually isomorphic, so $H^t(\hat X_Y,\eF)=H^t(Y_m,\eF)$ for $m$ large enough. Thus $H^t(X,\eF){\to} H^t(Y_m,\eF)$ are isomorphisms, so $H^t(\eI_Y^m\otimes \eF)=0$. 
It remains $t=p$. For $m\ges\cst{\eF}$, 
$H^p(\eI_Y^{m+1}\otimes\eF){\to} H^p(\eI_Y^{m}\otimes\eF)$ are injective, eventually isomorphic to a vector space $H(p,\eF)$. The previous step and the sequence $0{\to}\eI_Y^{m+1}{\to}\eO_X{\to}\eO_{Y_m}{\to}0$ imply that $H(p,\eF)=\Ker\big(H^p(X,\eF){\to} H^p(\hat X_Y,\eF)\big)=0$. 
\end{m-proof}

\begin{m-proposition}\label{prop:q1}
We consider the conditions: 
$$\text{{\rm(A)} $Y$ is a $(\dim Y-p)$-ample subscheme;
\qquad{\rm(P)} $Y$ is $p^\pos$.}$$
The following statements hold:
\begin{enumerate}[leftmargin=5ex]
\item[\rm(i)]
If $\tld X$ is a Cohen-Macaulay scheme, then $\,\text{\rm(A)}\,\Rightarrow\,\text{\rm(P)}$; 
\item[\rm(ii)]
If $\tld X$ is a Gorenstein, then $\,\text{\rm(A)}\,\Leftrightarrow\,\text{\rm(P)}$. 
\item[\rm(iii)] 
For $X$ smooth and $Y$ lci, one has the equivalence: 
\begin{equation}\label{eq:equiv-p}
\text{$Y$ is $p^\pos$}
\quad\Leftrightarrow\quad 
\bigg\{\begin{array}{l}
\text{the normal bundle $\eN_{Y/X}$ is $(\dim Y-p)$-ample,} 
\\[1ex] 
\text{the cohomological dimension $\cd(X\sm Y)\les\dim X-(p+1)$}. 
\end{array}
\end{equation}
\end{enumerate}
\end{m-proposition}

\begin{proof}
(i) Since $\eO_{\tld X}(-E_Y)$ is relatively ample, for $m\gg0$, it holds: 
$$
H^t(X,\eF\otimes\eI_Y^m)
\cong H^t(\tld X,\eF\otimes\eO_{\tld X}(-mE_Y))
\cong H^{\dim X-t}(\tld X,\omega_{\tld X}\otimes\eF^\vee\otimes\eO_{\tld X}(mE_Y)).
$$

\nit(ii) The equation above shows that~\eqref{eq:q11} holds for $\omega_{\tld X}\otimes\eL$, $\eL\in\Pic(X)$; let us prove for a coherent sheaf $\tld\eF$ on $\tld X$. Take $\eA\in\Pic(X)$ ample such that $\eA(-E_Y)$ is ample on $\tld X$, and $c>0$ such that $(\tld\eF\otimes\omega_{\tld X}^{-1})\otimes\eA(-E_Y)^c$ is globally generated. The recursion~\cite[Lemma~2.1]{ottm} applied to 
$\;0\to\tld\eF_1:=\Ker(\veps)\to
\big(\omega_{\tld X}\otimes\eA^{-c}\otimes\eO_{\tld X}(cE_Y)\big)^{\oplus N}
\srel{\veps}{\to}\tld\eF\to0,$ 
for some $N>0$, yields 
$H^j(\tld\eF(mE_Y))\subset H^{j+1}(\tld\eF_1(mE_Y))$, $j\ges\codim Y+q$, $m\gg0$. 

\nit(iii) Apply Proposition~\ref{prop:N+cd}, since $\tld X$ is Gorenstein.
\end{proof}

The Cohen-Macaulay (resp. Gorenstein) property of blow-ups has been investigated by several authors; combinatorial conditions are determined in~\cite{kaw,hyry}.  A situation which covers many geometric applications is when $X$ is Cohen-Macaulay (resp. smooth) and $Y$ is lci. 

The proposition above breaks the estimation of the amplitude into a local and a global problem. The former is easier but, in general, the cohomological dimension is difficult to control (cf.~\cite{andr+grau,ogus,falt-homog}); in~\cite{hlc+taj} we obtained upper bounds---those of interest---in the presence of affine stratifications. Below is a manageable situation which will be used in Section~\ref{sct:expl}.

\begin{m-proposition}\label{prop:x-b}
Let $V$ be a projective scheme, $\tld X\srel{\phi}{\to}V$ a morphism such that $\eO_{\tld X}(E_Y)$ is $\phi$-relatively ample. Then $Y\subset X$ is $q$-ample, for 
$\,q:=1+\dim\phi(\tld X)-\codim_X(Y).$
\end{m-proposition}

\begin{proof}
For a coherent sheaf $\tld\eF$ on $\tld X$, one has $R^t\phi_*(\tld\eF\otimes\eO_{\tld X}(mE_Y)){=}0,\,t>0, m\gg0$, so\\ 
$H^{j}\big(\tld X,\tld\eF\otimes\eO_{\tld X}(mE_Y)\big){=}H^j\big(\,V,\phi_*(\tld\eF\otimes\eO_{\tld X}(mE_Y))\,\big){=}0,\,j\ges q+\delta>\dim\phi(\tld X).$  
\end{proof}


\subsection{Elementary operations}\label{sct:N+cd}

We study the behaviour of partial ampleness under various natural operations: intersection, pull-back, product. 

\begin{m-proposition}\label{prop:pull-back}
Let $X'\srel{f}{\to}X$ be a morphism, $d:=$ the maximal dimension of its fibres, and $Y':=X'\times_XY=f^{-1}(Y)\subset X'$. 
\begin{enumerate}[leftmargin=5ex]
\item[\rm(i)] 
If $Y\subset X$ is $q$-ample, then $Y'\subset X'$ is $(q+d)$-ample. 
\item[\rm(ii)]
If $f$ is flat and surjective and $Y$ is $p^\pos$, then so is $Y'$. 
\end{enumerate}
\end{m-proposition}

\begin{m-proof}
(i) The universality property of the blow-up yields the commutative diagram 
$$
\xymatrix@R=1.2em@C=4em{
\tld X'=\Bl_{Y'}(X')\ar[d]\ar[r]^-{\tld f}&\tld X=\Bl_{Y}(X)\ar[d]
&
\\
X'\ar[r]^-f&X& \tld f^*\eO_{\tld X}(E_{Y})=\eO_{\tld X'}(E_{Y'}).
}
$$
Moreover, $\tld X'\subset X'\times_X\tld X$ is a closed subscheme, so the maximal dimension of the fibres of $\tld f$ is still $d$. For a coherent sheaf $\tld\eG$ on $\tld X'$, the projection formula implies  
$\,
R^i\tld f_*(\tld\eG\otimes\eO_{\tld X'}(mE_{Y'}))
=R^i\tld f_*\tld\eG\otimes\eO_{\tld X}(mE_{Y}),\; R^{j}\tld f_*\tld\eG=0,\,j>d,
$ 
and the condition~\eqref{eq:q11} follows. 

\nit(ii) Since $f$ is flat, the diagram is Cartesian. Take $\tld\eA'\in\Pic(\tld X')$ ample such that $R^j\tld f_*\tld{\eA'}^k=0,$ for $k,j>0$; then $\tld\eF_k:=\tld f_*\tld{\eA'}^k$ is locally free on $\tld X$, by Grauert's criterion. For $t\les p$, one has $H^t(\tld X',\tld{\eA'}^{k}(-mE_{Y'}))=H^t(\tld X,\tld\eF_k(-mE_Y))$, and we conclude by~\ref{prop:p>0}. 
\end{m-proof}

\begin{m-proposition}\label{prop:product}
Suppose $X$ is smooth. 
\begin{enumerate}[leftmargin=5ex]
\item[\rm(i)] 
Let $Y_1,Y_2\subset X$ be respectively $q_1$-, $q_2$-ample lci subvarieties such that 
$\;\codim(Y_1\cap Y_2)=\codim(Y_1)+\codim(Y_2).$ 
Then $Y_1\cap Y_2\subset X$ is $(q_1+q_2)$-ample. 
\item[\rm(ii)]
Suppose $Y_j\subset X_j$ are lci and $p_j^\pos$, for $j=1,2$. Then $Y_1\times Y_2\subset X_1\times X_2$ is $\min\{p_1,p_2\}^\pos$.
\end{enumerate}
\end{m-proposition}

\begin{m-proof}
(i) Note that $Y_1\cap Y_2$ is lci in $X$; the sub-additivity property~\cite[Theorem 3.1]{arap} applied to the normal bundle sequence implies that $\eN_{Y_1\cap Y_2/X}$ is $(q_1+q_2)$-ample. The Mayer-Vietoris sequence for $X\sm(Y_1\cup Y_2)$ yields the bound on the cohomological dimension.

\nit(ii) The inequality $\cd(X_1\times X_2\sm Y_1\times Y_2)<\dim(X_1\times X_2)-\min\{p_1,p_2\}$ follows from the Mayer-Vietoris sequence for $X_1\times X_2\sm Y_1\times Y_2=\big((X_1\sm Y_1)\times X_2\big)\cup\big(X_1\times(X_2\sm Y_2)\big)$. 

It remains to show that $\eN_{Y_1\times Y_2/X_1\times X_2}$ is $q$-ample, for $q=\dim Y_1+\dim Y_2-\min\{p_1,p_2\}$. We use the equivalent characterization in Definition~\ref{def:q-line}: let $\eA_1,\eA_2$ be ample line bundles on $X_1,X_2$, respectively, and $\eA_1\boxtimes\eA_2$ the tensor product of their pull-backs. Then, for $a\gg0$, $k\ges 1$, and  $t>\max\{q_1+\dim X_2,q_2+\dim X_1\}$, it holds:\\[1ex]
$\begin{array}{l}
H^t\big(Y_1\times Y_2,(\eA_1^{-k}\boxtimes\eA_2^{-k})\otimes\Sym^a(\eN_{Y_1/X_1}\boxplus\eN_{Y_2/X_2})\big)
\\[1.5ex]
\null\kern1em=
\uset{\genfrac{}{}{0pt}{}{t_1+t_2=t,}{a_1+a_2=a}}{\bigoplus}
H^{t_1}\big(Y_1,\eA_1^{-k}\otimes\Sym^{a_1}(\eN_{Y_1/X_1})\big)\otimes
H^{t_2}\big(Y_2,\eA_2^{-k}\otimes\Sym^{a_2}(\eN_{Y_2/X_2})\big)=0. 
\end{array}$
\end{m-proof}


\subsection{Weak positivity}\label{ssct:aprox-q}

Our goal is to prove a transitivity result for the $p^\pos$-property. 

\begin{m-definition}\label{def:ap>0} 
The subscheme $Y\subset X$ is $p^\apos$---that is, \emph{weakly $p$-positive}---if there is a decreasing sequence of sheaves of ideals $\{\eJ_m\}_{m}$ with the following properties:  
\begin{equation}\label{eq:ap}
\begin{array}{lll}
\bullet\;&
\forall\,m,n\ges1,\;\exists\,m'>m,\,n'>n\text{ such that }
\eJ_{m'}\subset\eI_Y^m,\;\eI_Y^{n'}\subset\eJ_{n};
&
\\[1ex] 
\bullet\;&
\text{for any locally free sheaf $\eF$ on $X,\;$}
\\[.5ex]&
\exists\,\cst{\eF}\ges1\text{ such that }H^t\big(X,\eF\otimes{\eJ_{m}}\big)=0,
\;\forall\,t\les p\;\forall\,m\ges\cst{\eF}.
&\kern10ex\null
\end{array}
\end{equation}
Obviously, $Y$ is $p^\apos$ if and only if so is $Y_{\text{red}}$. 
\end{m-definition}

\begin{m-lemma}\label{lm:cd-apos}
Suppose $X$ is smooth, $Y$ is lci and $p^\apos$. Then $\,\cd(X\sm Y){\les}\dim X{-}(p+1).$ 
\end{m-lemma}

\begin{m-proof}
Since $X$ is smooth, \cite[Proposition III.3.1]{hart-as} states: 
$$
\cd(X\sm Y)< c\;\Leftrightarrow\; 
H^t(X\sm Y,\eL)=0,\,\forall \eL\in\Pic(X),\,\forall t\ges c. 
$$
We have $X\sm Y\cong\tld X\sm E_Y$ and $H^t(\tld X\sm E_Y,\eL)=\varinjlim H^t(\tld X,\eL(mE_Y))$, cf.~\cite[(5.1)]{ottm}. Since  $\Pic(\tld X)\cong\Pic(X)\oplus\mbb ZE_Y$, one has $\omega_{\tld X}\otimes\eL^{-1}\cong \eM(lE_Y)$ for some $\eM\in\Pic(X)$, $l\in\mbb Z$. By applying Serre duality---as $\tld X$ is Gorenstein---we find:\\ 
$\varprojlim H^j(X,\eM\otimes\eI_Y^m)=\varprojlim H^j(X,\eM\otimes\eJ_n)=0, \,\forall\, \eM\in\Pic(X),\,j\les p.$ 
\end{m-proof}

\begin{m-corollary}\label{cor:pic}
Let $Y\subset X$ be complex, smooth varieties. If $Y$ is $p^\apos$, then the following statements hold: 
\begin{enumerate}[leftmargin=5ex]
\item[\rm(i)]  
$H^t(X;\bbQ)\to H^t(Y;\bbQ)\;\text{is}\;\biggl\{
\begin{array}{rl}
\text{an isomorphism, for}&t\les p-1;
\\[1ex] 
\text{injective, for}&t=p.
\end{array}$ 

\item[\rm(ii)] For $p\ges 3$, the following maps are isomorphisms: 
$$\Pic(X)\otimes\bbQ\to\Pic(\hat{X}_Y)\otimes\bbQ\to\Pic(Y)\otimes\bbQ.$$ 

\item[\rm(iii)] 
For $p=\dim Y-q$ and $X,Y$ as in~\ref{prop:x-b}, the previous properties hold with $\mbb Z$-coefficients. 
\end{enumerate}
\end{m-corollary}

\begin{m-proof}
(i) Use the previous lemma and~\cite[Corollary 5.2]{ottm}.

\nit(ii) One has $H^j(\eO_X)\srel{\cong}{\to}H^j(\hat X_Y;\cO_{\hat X_Y})$, $j=1,2$. On the other hand, the sequence 
$$
H^1(Y;\mbb Z)\to H^1(\hat X_Y;\cO_{\hat X_Y})\to \Pic(\hat X_Y)
\to H^2(Y;\mbb Z)\to H^2(\hat X_Y;\cO_{\hat X_Y})
$$ 
is exact (cf.~\cite[Lemma 8.3]{hart-cdav}). Now use (i) and the exponential sequences of $X,Y$.

\nit(iii) Note that $\eO_{\tld X}(E_Y)$ is $\dim\phi(\tld X)$-positive, so~\ref{thm:mats-bott}(iii) applies. Indeed, consider an embedding $\tld X\srel{\iota}{\to}\mbb P^N\times V$ over $V$, such that $\eO_{\tld X}(m_0E_Y)=\iota^*(\eO_{\mbb P^N}(1)\boxtimes\euf M)$, for some $m_0>0$, $\euf M\in\Pic(V)$. Now take the Fubini-Study metric on $\eO_{\mbb P^N}(1)$ and an arbitrary on $\euf M$. 
\end{m-proof}

\begin{m-proposition}\label{prop:approx-p}
Let $X$ be smooth. Suppose $Z\subset Y$ is $p^\pos$, $Y\subset X$ is $r^\pos$, and they are both irreducible and lci. Then the following statements hold:
\begin{enumerate}[leftmargin=5ex]
\item[\rm(i)] 
$Z\subset X$ is $\big(p-(\dim Y-r)\big)^\pos$; more precisely, one has:
\begin{equation}\label{eq:trans-p}
\bigg\{
\begin{array}{l}
\text{$\eN_{Z/X}$ is $\big(\dim Y+\dim Z-(r+p)\big)$-ample,}
\\[1ex] 
\text{$\cd(X\sm Z)\les\dim X-(\min\{r,p\}+1)$}.
\end{array}
\end{equation}

\item[\rm(ii)]
If $Z, Y$ are smooth, then $Z\subset X$ is ${\mn{\{p,r\}}}^\apos$.
\end{enumerate}
\end{m-proposition}

\begin{m-proof}
(i) The first claim follows from~\ref{prop:N+cd}. For the second, let $U_Z{:=}X\sm Z,\, U_Y{:=}X\sm Y$, and $\eG$ be  a coherent sheaf on $X$. The left- and right-hand side of the exact sequence 
$$
\ldots\to H^i_{Y\sm Z}(U_Z,\eG)\to H^i(U_Z,\eG)\to H^i(U_Y,\eG)\to\ldots,
$$
vanish for $i\ges\dim X-p$ and $i\ges \dim X-r$, respectively (cf.~\cite[Proposition 6.4]{ottm}). 

\nit(ii) Both $Y,Z$ are ${\mn{\{p,r\}}}^\pos$, so we may assume without loss of generality that $p=r$. Consider $\xi_1,\dots,\xi_u,\zeta_1,\dots,\zeta_v\in\eO_{X,z}$, whose images in $\hat{\cal O}_{X,z}$ yield independent variables, such that $\eI_{Y,z}=\lran{\bsymb{\xi}}=\lran{\xi_1,\dots,\xi_u}$ and  $\eI_{Z,z}=\lran{\bsymb{\xi},\bsymb{\zeta}}=\lran{\xi_1,\dots,\xi_u,\zeta_1,\dots,\zeta_v}$. For $l\ges a$, a direct computation yields 
$\,\eI_{Y,z}^a\cap\eI_{Z,z}^l
=\ouset{i=a}{l}{\sum}\lran{\bsymb{\xi}}^{i}\cdot\lran{\bsymb{\zeta}}^{l-i}
=\eI_{Y,z}^a\cdot\eI_{Z,z}^{l-a},$
which implies  $(\eI_{Z,z}^l+\eI_{Y,z}^{a})/{\eI_{Y,z}^{a}}\cong\eI_{Z,z}^{l}/\eI_{Y,z}^a\cdot\eI_{Z,z}^{l-a}$.
We obtain the exact sequences: 
\begin{equation}\label{eq:la}
0\to
\frac{\eI_Y^a}{\eI_Y^{a+1}}\otimes\biggl(\frac{\eI_Z}{\eI_Y}\biggr)^{l-a}
\!\to
\frac{\eI_Z^l+\eI_Y^{a+1}}{\eI_Y^{a+1}}
\to
\frac{\eI_Z^l+\eI_Y^{a}}{\eI_Y^{a}}
\to0,\quad\forall\,l\ges a+1.
\end{equation}
The left side is an $\eO_Y$-module: $\eI_Z/\eI_Y=\eI_{Z\subset Y}$ is the ideal of $Z\subset Y$; $\;\eI_Y^a/\eI_Y^{a+1}=\Sym^a\eN_{Y/X}^\vee$. 

Let $\eF$ be locally free on $X$. By the $p^\pos$-property, there is a linear function $l(k)=\cst{}_1\!\cdot\,k+\cst{}_2$ (with $\cst{}_1,\cst{}_2$ independent of $\eF$) and $k_\eF,l_\eF\in\mbb N$, such that: 
$$
\begin{array}{rl}
H^t(\eF\otimes\eI_{Y}^k)=0,
&
\;\forall\,t\les p,\;\forall\,k\ges k_\eF,
\\[1ex]
H^t(\eF_Y\otimes\eI_{Z\subset Y}^{l})=0,
&
\;\forall\,t\les p,\;\forall\,l\ges l_\eF,
\\[1ex] 
H^t(\eF_Y\otimes\Sym^a\eN_{Y/X}^\vee\otimes\eI_{Z\subset Y}^{l-a})=0,
&
\;\forall\,t\les p,\;\forall\,a\les k,\;\forall\,l\ges l(k).
\end{array}
$$
The last claim is a consequence of the uniform $q$-ampleness and the sub-additivity property of the amplitude (cf.~\cite[Theorems 7.1]{tot}, \cite[Theorem 3.1]{arap}): 
\begin{itemize}[leftmargin=5ex]
\item[--]
There is a function $\text{linear}(r)$ such that, for any locally free sheaf $\eF$ whose regularity satisfies $\max\{1, \reg(\eF_Y)\}\les r$, it holds: 
$\;H^t(\eF_Y\otimes\eI_{Z\subset Y}^l)=0,\;\forall\,t\les p,\; l\ges{\rm linear}(r).$
\item[--] 
If $a\les k$, then ${\rm reg}(\eF_Y\otimes\Sym^a\eN_{Y/X}^\vee)\les{\rm linear}(k)$. 
\end{itemize}
Recursively for $a=1,\ldots,k$, and starting by $\frac{\eI_Z^l+\eI_Y}{\eI_Y}=\eI_{Z\subset Y}^l$, \eqref{eq:la} yields: 
\\[.5ex] \centerline{
$
H^t\Big(
\eF\otimes\frac{\eI_Z^l+\eI_Y^{k}}{\eI_Y^{k}}
\Big)=0,\;\forall t\les p,\;\forall\,l\ges l(k). 
$
}\\[.5ex]
Now tensor $0\to\eI_Y^k\to\eI_Z^l+\eI_Y^{k}\to\frac{\eI_Z^l+\eI_Y^{k}}{\eI_Y^{k}}\to0$ by $\eF$ and deduce: 
$$
H^t\big(\eF\otimes(\eI_Y^k+\eI_Z^l)\big)=0,
\;\forall\,t\les p,\;\;\forall\,k\ges k_\eF,\;\forall\,l\ges l(k).
$$
The subschemes defined by $\eI_Y^k+\eI_Z^l$ are `asymmetric' thickenings of $Z$ in $X$. The ideals $\eJ_k:=\eI_Y^k+\eI_Z^{k+l(k)}$ satisfy~\eqref{eq:ap}: 
$\eJ_{k'}\subset\eI_Z^k,\;k'\ges k,\;\; \eI_Z^{m'}\subset\eJ_m,\;m'\ges m+l(m).$
\end{m-proof}


\section{Connectedness properties}\label{ssct:G3}

\begin{m-notation}\label{not:VfX}
Let $V,X$ be irreducible projective varieties, $V\srel{f}{\to} X$ a morphism, and $Y\subset X$ a closed subscheme. 
\end{m-notation}

The issue regarding the connectedness of pre-images of subschemes by morphisms was raised by Fulton-Hansen in the late~70s. 

\begin{conj-nono}{(cf.~\cite[p.\,161]{fult+hans})} 
Suppose that $\dim f(V)+\dim Y>\dim X$ and the normal bundle $\eN_{Y/X}$ is ample. Then $f^{-1}(Y)$ is connected. 
\end{conj-nono}

Despite its elementary nature, it turns out that the question is surprisingly difficult to answer. It is known that, in this form, the conjecture is false; a counterexample can be found in \cite{hart-as}. However, it does hold for subvarieties of various homogeneous spaces (cf.~\cite{fult+hans,hans,falt-homog,debr2}), so it is interesting to find a framework which yields a positive answer.


\subsection{Connectedness of pre-images}\label{ssct:conn}

\begin{m-theorem}\label{thm:f-1}
Suppose $\cd(X\sm Y)\leq\dim(X)-2$ and $f$ is surjective. Then $f^{-1}(Y)$ is connected, in particular so is $Y$.
\end{m-theorem}
The statement generalizes~\cite[Corollary III.3.9]{hart-as} in two directions. First, it allows morphisms into the picture; this is important, taking into account the Fulton-Hansen-problem. Second, there is no assumption on the smoothness of the varieties; this is crucial, since one can not control the regularity of the image of an arbitrary morphism.

\begin{m-proof}
We may assume that $V$ is smooth. Otherwise, let $V'\srel{\si}{\to}V$ be a (surjective) resolution of singularities; if $(f\si)^{-1}(Y)$ is connected, then so is $f^{-1}(Y)=\si\big((f\si)^{-1}(Y)\big)$. 

In order to prove that $Z:=f^{-1}(Y)$ is connected, it suffices to show that 
\begin{equation}\label{eq:res}
\res_Z:H^0(V,\eO_V)\to H^0(\hat V_{Z},\eO_{\hat V_Z})
\end{equation}
is an isomorphism. By formal duality~\cite[Theorem~III.3.3]{hart-as}, the right-hand side is isomorphic to $H_Z^{\dim V}(V,\omega_V)^\vee$ and the dual of $\res_Z$ fits into the exact sequence: 
$$
H^{\dim V-1}(V\setminus Z,\omega_V)\to H_Z^{\dim V}(V,\omega_V)\to H^{\dim V}(V,\omega_V)\to H^{\dim V}(V\setminus Z,\omega_V).
$$
The rightmost cohomology group vanishes, by Lichtenbaum's theorem. 

We claim that the leftmost group vanishes too. This a consequence of Leray's spectral sequence for $f$, combined with Koll\'ar's higher direct image theorem~\cite[Theorem 2.1]{kolr-I}. Indeed, the cohomology group $H^{a+b}(V\setminus Z,\omega_V)$ can be computed using the spectral sequence whose $E_2$-term is $H^a(X\setminus Y,R^bf_*\omega_V)$. With the \textit{ad hoc} notation 
$$
v:=\dim V,\; x:=\dim X,\;\text{so}\;v-x=\dim(\text{generic fibre of}\;f),
$$
Koll\'ar's theorem states that $R^bf_*\omega_V=0,\;\text{for}\;b\ges v-x+1.$ 

The restriction of $f$ to $V\setminus Z$ is proper, so the higher direct images are coherent (in fact torsion free, by \lcit). The assumption on the cohomological dimension of $X\setminus Y$ yields $H^a(X\setminus Y,R^bf_*\omega_V)=0,\;\text{for}\;a\geq x-1.$
\end{m-proof}

\begin{m-corollary}\label{cor:etale}
The induced homomorphism $\pi_1^{alg}(f^{-1}(Y))\to\pi_1^{alg}(V)$ between the algebraic fundamental groups is surjective.
\end{m-corollary}

\begin{m-proof}
For any \'etale morphism $W\srel{g}{\to}V$, $(fg)^{-1}(Y)=g^{-1}\big(f^{-1}(Y)\big)$ is connected. 
\end{m-proof}

In general, for arbitrary $V$, there is no control on the homomorphism~\eqref{eq:res}. 

\begin{m-proposition}\label{thm:rtl}
Let the situation be as above and suppose moreover that $V$ is normal and has rational singularities. Then $\res_Z:H^0(V,\eO_V)\to H^0(\hat V_{Z},\eO_{\hat V_Z})$ is an isomorphism. 
\end{m-proposition}
Subschemes satisfying this property are called G1 in Hironaka-Matsumura~\cite{hir+mats}.\smallskip 

\begin{m-proof}
The argument is the same as above. Kempf's criterion implies that $V$ is Cohen-Macaulay and, at the first step of the previous proof, the resolution $V'\srel{\si}{\to} V$ has the property that $\si_*\omega_{V'}=\omega_V$. Consequently, formal duality---needed for dualizing~\eqref{eq:res}---holds on $V$.  (Cohen-Macaulayness suffices for the Serre duality in the proof of~\cite[Theorem~III.3.3]{hart-as}.) For $Z':=\si^{-1}(Z)$, one has: 
$\,H^{j}(V'\setminus Z',\omega_{V'})\cong H^{j}(V\setminus Z,\omega_V),\;j=v-1,v;$
this transfers the computation from $V$ to $V'$, which is smooth.
\end{m-proof}

\begin{m-remark}\label{rmk:optim}
\begin{enumerate}[leftmargin=5ex]
\item 
Lichtenbaum's theorem states that $\cd(X\sm Y)\leq x-1$, but it is unclear when is maximal. Our result implies that, contrary to the intuition, the cohomological dimension does not drop by removing effective divisors, the reason is not due to the existence of `disjoint divisors'. Indeed, suppose $\cd(X\sm Y)=x-1$ and $D\subset X$ is a (complete) effective divisor, disjoint of $Y$. Then it still holds $\cd(X\sm(Y\cup D))=x-1$; otherwise $Y\cup D$ would be connected, contradicting the hypothesis.

The observation is false if the divisor is allowed to intersect $Y$: $\cd(\mathbb P^2\sm[1{:}0{:}0])=1$ and it remains the same by removing a line disjoint of $[1{:}0{:}0]$. However, by removing a line passing through the point, one obtains the affine $2$-plane, whose cohomological dimension vanishes.
\item 
The Fulton-Hansen-conjecture is false in general; a counterexample is due to Hartshorne (cf.~\cite[pp.\,199]{hart-as}). Let $V{\srel{f}{\to}}X$ be an \'etale (surjective) morphism, $Y'$ a $\delta$-co\-dimensional, general complete intersection in $V$, with $2\delta\,{>}\dim V$; let $Y{:=}f(Y')$. Then $\eN_{Y/X}$ is ample, $f^{-1}(Y')$ is disconnected. What (necessarily) fails is $\cd(X\sm Y){\leq}\dim X{-}2$. This shows that Theorem~\ref{thm:f-1} is optimal.
\end{enumerate}
\end{m-remark}


\subsection{Application to partially ample subvarieties}\label{ssct:applic}

We start discussing the equidimensionality of partially ample subvarieties. In general, they are not so: consider $X:=\mbb P^2$ and $Y:=\{x=0\}\cup\{y=z=0\}$. Then $\tld X=\Bl_Y(\mbb P^2)$ is isomorphic to the blow-up $\wtld{\mbb P^2}$ of $\mbb P^2$ at $[1:0:0]$, with exceptional divisor $E$, and $\eO_{\tld X}(E_Y)=\eO_{\mbb P^2}(1)\otimes\eO_{\wtld{\mbb P^2}}(E)$. A short computation shows that $Y$ is $1$-ample. 

\begin{m-lemma}\label{lm:conn}
Let $X$ be a projective variety, and suppose that $X, Y$ are both Cohen-Macaulay. If $Y$ is either $1^\apos$ or $(\dim Y-1)$-ample in $X$, then $Y$ is equidimensional.
\end{m-lemma}

\begin{m-proof}
Suppose $Y$ is $1^\apos$. Then it is G1: the isomorphism $H^0(\eO_X)\cong H^0(\eO_{\hat X_Y})$ holds `in finite time': with notation~\ref{def:ap>0}, one has $H^0(\eJ_m)=H^1(\eJ_m)=0$, $m\gg0$. In particular, $Y$ is connected. 
The same holds if $Y$ is $(\dim Y-1)$-ample. Now conclude by using the unmixedness property, that local Cohen-Macaulay rings are equidimensional.
\end{m-proof}

Recall that partial ampleness yields an upper bound for the cohomological dimension of the complement of a subvariety. Thus we obtain a convenient class of subvarieties for which Theorem~\ref{thm:f-1} does apply. Below is our main result: partial ampleness is a \emph{numerical condition} which ensures that pre-images are connected. We stress that, in the generality below, there are \emph{no similar statements} in the literature.

\begin{m-theorem}\label{thm:q}
Let the situation be as in~\ref{not:VfX}.
\begin{enumerate}[leftmargin=5ex]
\item[\rm(i)] 
Let $Y\subset X$ be a $\big(\dim f(V)+\dim(Y)-\dim(X)-1\big)$-ample closed subscheme. 
Then the following properties hold:
$$\text{$f^{-1}(Y)\subset V$ is connected;\qquad $\pi_1^{alg}\big(f^{-1}(Y)\big)\to\pi_1^{alg}(V)$ is surjective.}$$
\item[\rm(ii)] 
In particular, suppose $Y\subset X$ is a $(\dim Y-1)$-ample subscheme and $f$ is surjective. Then $f^{-1}(Y)$ is connected. 
\end{enumerate}
\end{m-theorem}

\begin{m-proof}
Observe that $\cd\big(\,f(V)\sm(Y\cap f(V))\,\big)\les\cd(X\sm Y)\les\dim f(V)-2$, and apply~\ref{thm:f-1} to the surjective morphism $f:V\to f(V)$.
\end{m-proof}

\begin{m-corollary}\label{cor:fh}
Let $V, Y$ be closed subschemes of $X$. Suppose $V$ is connected and $Y$ is $q$-ample in $X$, with $0\les q\les\dim V+\dim Y-\dim X-1.$ 

Then $Y\cap V$ is non-empty and $\big(\dim(Y\cap V)-1\big)$-ample in $V$, hence connected. 
\end{m-corollary}
The statement reminds a problem of Hartshorne \cite[Ch.~III, Conjecture~4.5]{hart-as}, concerning the connectedness and non-emptiness of the intersection of smooth subvarieties with ample normal bundles. According to the survey~\cite{petr}, this issue is currently still wide open.\smallskip 

\begin{m-proof}
The intersection  $Y\cap V$ is non-empty, because 
$$\;\cd(V\sm Y\cap V)\les\cd(X\sm Y)\les\codim_XY+q-1=\dim V-2,\quad\text{(cf.~\ref{prop:N+cd}).}$$
Actually one has $\dim(Y\cap V)\ges1$, since otherwise $\cd(V\sm Y\cap V)=\dim V-1$: there are effective divisors avoiding a finite number of points. 
\end{m-proof}

\begin{m-example}
Suppose $Y\subset X$ is ample. Then, for any subvariety $V\subset X$ of dimension at least $\codim(Y)+1$, the intersection $Y\cap V$ is non-empty---this is already proved in~\cite{ottm}---and also connected. The Lefschetz-type hyperplane theorem in \ocit, Corollary~5.2---in particular the connectedness of the intersection---requires the smoothness of $V\sm(Y\cap V)$. 

Therefore, by pursuing this path, one can not deduce the connectedness of pre-images, as we do in Theorem~\ref{thm:q}, because it is not possible to control the smoothness of the image of an arbitrary morphism.
\end{m-example}

So far we used only the bound on the cohomological dimension of $X\sm Y$. The partial ampleness of $Y\subset X$ actually carries more information.

\begin{m-theorem}
\begin{enumerate}[leftmargin=5ex]
\item[\rm(i)] 
Let the situation as in~\ref{cor:fh}, $X,V$ smooth; suppose $Y{\subset} X$ and $Y{\cap} V{\subset} V$ are lci. (The intersection is automatically lci if  $\codim_X(V{\cap} Y){=}\codim_XV{+}\codim_XY.$) 
\item[] 
Then $V{\cap} Y$ is G3 in $V$. In particular, a $1^\pos$ lci subscheme of a smooth variety is G3.
\item[\rm(ii)] 
Let the situation as in~\ref{not:VfX}, with $f(V)$ is smooth and $Y\cap f(V)$ is lci. 
\item[] 
If $Y$ is $\big(\dim f(V)+\dim(Y)-\dim(X)-1\big)$-ample, then $f^{-1}(Y)$ is G3 in $V$.
\end{enumerate}
\end{m-theorem}
For the definition of the G3-property, the reader in invited to consult~\cite{hir+mats}.\smallskip

\begin{m-proof}
(i) First note that if $V\cap Y$ has the expected codimension, then each of its components has at most that codimension, so $V\cap Y$ is equidimensional; thus $V\cap Y\subset V$ is lci. 

Back to the general case, the commutative diagram below shows that the exceptional divisor $E_{Y\cap V}$ is $(\dim V-2)$-ample:
$$
\xymatrix@R=1.2em{
\Bl_{Y\cap V}(V)\;\ar@{^(->}[r]\ar[d]&\Bl_{Y}(X)\ar[d]
\\ 
V\;\ar@{^(->}[r]&X.
}
$$
Hence $\eN_{V{\cap} Y/V}$ is $\big(\dim(V{\cap} Y)-1\big)$-ample by~\ref{prop:q1}(iii), so $V{\cap} Y$ is G2 in $V$ (cf.~\cite[\S3]{hlc-subvar}). As $\cd(V\sm V{\cap} Y)\les\dim V-2$, Speiser's result~\cite[Corollary V.2.2]{hart-as} yields the conclusion. 

\nit(ii) The proof of~\ref{thm:q} above and the previous step imply that $Y\cap f(V)$ is G3 in $f(V)$. It remains to apply~\cite[Theorem~2.7]{hir+mats}, since $f^{-1}(Y)=f^{-1}(Y\cap f(V))$. 
\end{m-proof}
In Hartshorne's counterexample~\ref{rmk:optim}, $Y{\subset} X$ is not G3, but it is G2. Also, $Y$ does not possess the $1^\pos$-property. This indicates that, for being G3, the $1^\pos$-property is close to optimal; see~\cite[Proposition~V.2.1]{hart-as}. It has the advantage to be a numerical condition.


\section{Examples of partially ample subvarieties}\label{sct:expl}

We show that partially ample subvarieties occur in a variety of situations:
\begin{enumerate}
\item zero loci of sections in vector bundles;
\item sources of Bialynicki-Birula decompositions;
\item subvarieties of rational homogeneous varieties.
\end{enumerate}


\subsection{Vanishing loci of sections}\label{ssct:glob-gen}

Throughout this section, $\eN$ is a vector bundle of rank $\nu$ on the smooth projective variety $X$.

\subsubsection{$q$-ample vector bundles}\label{sssct:q-vb}

\begin{m-proposition} \label{prop:q21} 
Suppose $\eN$ is $q$-ample and $Y$ is the zero locus of a \emph{regular} section in it. Then $Y\subset X$ is a $q$-ample subvariety.
\end{m-proposition}

\begin{proof} 
We verify~\eqref{eq:q12} for a vector bundle $\eF$ on $X$. Since $s$ is regular, $Y$ is lci, $\codim_X(Y)=\nu$, so~\ref{prop:q1}(ii) applies. One has the resolution (cf. \cite[Theorem 3.1]{bu+ei})
\begin{equation}\label{eq:koszul-m}
0\to L^\nu_m(\eN^\vee)\to\ldots\to L^j_m(\eN^\vee)\to\ldots\to 
\Sym^m(\eN^\vee)\srel{s^m\ort}{-\kern-1ex-\kern-1ex\lar}\eI_Y^m\to 0,\;\;\forall\,m\ges 1,
\end{equation}
where 
$L^j_m(\eN^\vee):=
\Img\Bigl(
\Sym^{m-1}(\eN^\vee)\otimes\overset{j}{\hbox{$\bigwedge$}}\,\eN^\vee 
\srel{\phi^j_m}{-\kern-1ex\lar} 
\Sym^{m}(\eN^\vee)\otimes\overset{j-1}{\hbox{$\bigwedge$}}\,\eN^\vee 
\Bigr), 1\les j\les\nu.$ 
The general linear group is linearly reductive and $\phi^j_m$ is equivariant, so $L^j_m(\eN^\vee)$ is a direct summand of $\Sym^m(\eN^\vee)\otimes\overset{j-1}{\hbox{$\bigwedge$}}\,\eN^\vee$. For $m\gg0$, one has: 
\\ \centerline{
$H^{t+j-1}(X,\eF\otimes\oset{j-1}{\bigwedge}\eN^\vee\otimes\Sym^m\eN^\vee)=0$,\; 
for $1\les j\les\nu,\;\;t+\nu-1\les\dim X-q-1$.
}\\[1ex] 
It follows that $H^t(X,\eF\otimes\eI_Y^m)=0$, for $0\les t\les\dim Y-q$.
\end{proof}


\subsubsection{Globally generated vector bundles}\label{ssct:x-b}

Henceforth we assume that $\eN$ is globally generated; thus the notions of $q$-ampleness and Sommese-$q$-ampleness agree (cf.~\ref{prop:q2}). 

Let $Y\subset X$ be lci of codimension $\delta$, the zero locus of $s\in\Gamma(\eN):=H^0(X,\eN)$. We \emph{do not require} $s$ to be regular, so we allow $\delta<\nu$. We are going to use~\ref{prop:x-b} to estimate the ampleness of $Y$. We observe that the blow-up fits into the diagram 
\begin{equation}\label{eq:tld-x}
\xymatrix@R=1.2em@C=1.75em{
\tld X\ar@{^(->}[r]\ar[d]_-\si\ar@<-5pt>[rrd]_-\phi
&
\mbb P(\eN)
{=}\,
\mbb P\Bigl(
\mbox{$\overset{\nu-1}{\bigwedge}$}\eN^\vee\otimes\det(\eN)
\Bigr)\ar@{^(->}[r]
&
X\times\mbb P\Bigl(
\mbox{$\overset{\nu-1}{\bigwedge}$}\Gamma(\eN)^\vee
\Bigr)\ar[d]
\\ 
X&&\mbb P:=\mbb P\Bigl(
\mbox{$\overset{\nu-1}{\bigwedge}$}\Gamma(\eN)^\vee
\Bigr),
}
\end{equation}
and it holds 
\begin{equation}\label{eq:o1}
\eO_{\tld X}(E_Y)=\eO_{\mbb P(\eN)}(-1)\big|_{\tld X}=
\bigl(\det(\eN)\boxtimes\eO_{\mbb P}(-1)\bigr)\big|_{\tld X}.
\end{equation}

\begin{m-proposition}\label{prop:p}
Suppose $\det(\eN)$ is ample. If the dimension of the generic fibre of $\phi$ over its image is $p+1$, then $\eO_{\tld X}(E_Y)$ is $\dim\phi(\tld X)$-positive, and $Y$ is $(\dim Y-p)$-ample. 
\end{m-proposition}

\begin{proof}
The assumptions of \ref{prop:x-b} are satisfied. 
\end{proof}

Note that the proposition applies also when only some symmetric power $\Sym^a\eN$ is globally generated. Then $s\in\Gamma(\eN)$ induces $s^a\in\Gamma(\Sym^a\eN)$ and $\eI_{\{s^a=0\}}=\eI_{\{s=0\}}^a$. By~\ref{prop:N+cd}(i), the amplitude of $\{s^a=0\}$ coincides with that of $\{s=0\}$.


\subsubsection{Special Schubert subvarieties of the Grassmannian}
\label{sssct:spec-grass}

Let $W\subseteq\Gamma(\eN)$ be a vector subspace generating $\eN$, $\dim W=\nu+u+1$. It is equivalent to a morphism $f:X\to\Grs(W;\nu)$ to the Grassmannian of $\nu$-dimensional quotients; $\det(\eN)$ is ample when $\vphi$ is finite. 

Henceforth let $X=\Grs(W;\nu)$; it is isomorphic to $\Grs(u+1;W)$, the $(u+1)$-dimensional subspaces of $W$; let $\eN$ be the universal quotient. The morphism $\phi$ in \eqref{eq:tld-x} is explicit: 
\begin{equation}\label{eq:q}
\mbb P(\eN)\to\mbb P,\quad 
(x,\lran{e_x})\mt \det(\eN_x/\lran{e_x})^\vee\subset
\oset{\nu-1}{\bigwedge}\eN_x^\vee\subset
\oset{\nu-1}{\bigwedge}W^\vee.
\end{equation}
($\lran{e_x}$ stands for the line generated by $e_x\in\eN_x$, $x\in\Grs(W;\nu)$.) 
The restriction to the Grassmannian corresponds to the commutative diagram  
\begin{equation}\label{eq:wn}
\xymatrix@R=1.2em{
0\ar[r]&\eO_{\Grs(W;\nu)}\ar[r]^-{s}\ar@{=}[d]&
W\otimes\eO_{\Grs(W;\nu)}\ar[r]\ar@{->>}[d]^-{\;\beta}&
W/\lran{s}\otimes\eO_{\Grs(W;\nu)}\ar[r]\ar@{->>}[d]&0
\\ 
&\eO_{\Grs(W;\nu)}\ar[r]^-{\beta s}&\eN\ar[r]&\eN/\lran{\beta s}\ar[r]&0.
}
\end{equation}
Thus $\phi$ is the desingularization of the rational map 
\begin{equation}\label{eq:q-grs}
g_s:\Grs(W;\nu)\dashto\Grs(W/\lran{s};\nu-1),\quad 
[W\surj N]\mt [W/\lran{s}\;\surj N/\lran{\beta s}],
\end{equation}
followed by the Pl\"ucker embedding of $\Grs(W/\lran{s};\nu-1)$. The indeterminacy locus of $\phi$ is $\Grs(W/\lran{s};\nu)$, so the latter is $u^\pos$ in $\Grs(W;\nu)$. The observation can be generalized. 
\begin{m-corollary}\label{cor:Yl}
For $\ell\les\nu$, fix an $\ell$-dimensional subspace $\Lambda_\ell\subset W$. Consider the Schubert subvariety 
$Y_\ell:=\{U\in\Grs(u+1;W)\mid U\cap\Lambda_\ell\neq0\}.$ 
Then $Y_\ell$ is  $\big(\ell(u+1)-1\big)^\pos.$ 
\end{m-corollary}
Thus the Chow ring of the Grassmannian is generated by partially ample subvarieties. 
\begin{m-proof}
Note that $Y_\ell$ is $(\nu-\ell+1)$-codimensional and is the vanishing locus of 
$$
s_\ell:\eO\cong\det(\Lambda_\ell\otimes\eO)\to\oset{\ell}{\bigwedge} W\otimes\eO\to\oset{\ell}{\bigwedge}\eN.
$$ 
We are in the situation~\ref{prop:x-b}. The diagram \eqref{eq:tld-x} corresponds to the rational map 
$$
\phi:\Grs(u+1;W)\dashto\Grs(u+1;W/\Lambda_\ell),\quad U\mt(U+\Lambda_\ell)/\Lambda_\ell, 
$$
followed by a large Pl\"ucker embedding; its indeterminacy locus is precisely $Y_\ell$. Since $\phi$ is surjective, a dimension counting yields the conclusion.
\end{m-proof}

\begin{m-remark}\label{rmk:sommese-weak}
\begin{enumerate}[leftmargin=5ex]
\item[\rm(i)] 
Propositions~ \ref{prop:q21} and~\ref{prop:p} deal with complementary situations: $\eO_{\mbb P(\eN^\vee)}(1)$ is the pull-back of an ample line bundle, while $\eO_{\tld X}(E_Y)$ is relatively ample.

\item[\rm(ii)] 
The criterion \ref{prop:q21} is not optimal: by~\ref{prop:q2}, for $X=\Grs(\nu+u+1;\nu)$, the universal quotient $\eN$ is $q$-ample, with $q=\dim\mbb P(\eN^\vee)-\mbb P^{\nu+u}=\dim X-(u+1)$. So $Y=\Grs(\nu+u;\nu)$, the zero locus of a section of $\eN$, is $(u+1-\nu)^\pos$; this may be negative and the estimate irrelevant. 

On the other hand, \ref{prop:x-b} implies that $Y$ is $u^\pos$. Moreover, for $\ell\neq1,\nu$, the section $s_\ell$ above is \emph{not regular}, so~\ref{prop:q21} does not apply, anyway.

\item[\rm(iii)] 
Subvarieties obtained as zero loci of sections in globally generated vector bundles and pull-backs of Schubert cycles appear in the recent work~\cite[\S3]{ful+leh}, in the definition of the pliant cone of a projective variety, which is a full-dimensional subcone of the nef cone---an object of central interest. The discussion above, together with Proposition~\ref{prop:pull-back}, implies that these elements of the pliant cone are in fact partially ample. 
\end{enumerate}
\end{m-remark}


\subsection{Sources of torus actions}\label{ssct:fixed} 

Let $X$ be a smooth projective variety with a faithful action $\lda:G_m\times X\to X$ of the multiplicative group $G_m=\kk^\times$. This determines the well-known Bialynicki-Birula---BB for short---decomposition of $X$ (cf.~\cite{bb}): 
\begin{itemize}[leftmargin=5ex]
	\item[$\bullet$] 
		The fixed locus $X^\lda$ of the action is a disjoint union 
		$\underset{s\in S_\BB}{\bigsqcup}\kern-1exY_s$ of smooth subvarieties. 
		For $s\in S_\BB$, 
		$Y_s^+:=\{x\in X\mid\underset{t\to 0}{\lim}\,\lda(t, x)\in Y_s\}$ 
		is locally closed in $X$ (a BB-cell) and it holds:
		$\;X=\underset{s\in S_\BB}{\bigsqcup}\kern-1exY_s^+.$
	\item[$\bullet$] 
		The \emph{source} $Y:=Y_{\rm source}$ and 
		the \emph{sink} $Y_{\text{\rm sink}}$ of the action are uniquely characterized 
		by the conditions: $Y^+=Y_{\rm source}^+\subset X$ is open 
		and $Y_{\rm sink}^+=Y_{\rm sink}$. 
\end{itemize}
A linearization of the action in a sufficiently ample line bundle yields a $G_m$-equivariant embedding $X\subset\mbb P^N_\kk$. There are homogeneous coordinates $\bz_{0}\in\kk^{N_0+1},\ldots,\bz_{r}\in\kk^{N_r+1}$ such that the $G_m$-action on $\mbb P^N_\kk$ is: 
\begin{equation}\label{eq:c*}
\lda\big(t,[\bz_{0},\bz_{1},\ldots,\bz_{r}]\big)
=[\bz_{0},t^{m_1}\bz_{1},\ldots,t^{m_r}\bz_{r}],\quad\text{with}\;0<m_1<\ldots<m_r. 
\end{equation}
The source and sink of $\mbb P^N, X$ are respectively: 
\begin{equation}\label{eq:YP}
\begin{array}{lcll}
\mbb P^N_{\rm source}=\{[\bz_{0},0,\ldots,0]\}, 
&&
\mbb P^N_{\rm sink}=\{[0,\ldots,0,\bz_{r}]\},
\\[1ex] 
Y=Y_{\rm source}=X\cap \mbb P^N_{\rm source},
&& 
Y_{\rm sink}=X\cap \mbb P^N_{\rm sink},
\\[1ex]
Y^+=X\cap (\mbb P^N_{\rm source})^+,
&&
(\mbb P^N_{\rm source})^+
=\{[\bz]=[\bz_{0},\bz_{1},\ldots,\bz_{r}]\mid\bz_{0}\neq 0\}.
\end{array}
\end{equation}

Let $m$ be the lowest common multiple of $\{m_\rho\}_{\rho=1,\dots,r}$ and $l_\rho:=m/m_\rho$. Let us denote $\bz_{\rho}^{l_\rho}:=(z_{\rho 0}^{l_\rho},\dots,z_{\rho N_\rho}^{l_\rho})$ and $\eI\subset\eO_{\mbb P^N}$ the sheaf of ideals generated by $\bz_1^{l_1},\dots,\bz_r^{l_r}$. The choice of the exponents makes the assignment 
\begin{equation}\label{eq:phi}
\mbb P^N\dashto\mbb P^{N'},\quad
[\bz_0,\bz_1,\dots,\bz_r]\mt[\bz_1^{l_1},\dots,\bz_r^{l_r}],
\end{equation}
well-defined. It defines a $G_m$-invariant rational map whose indeterminacy locus is the subscheme determined by $\eI$. Then $\eJ:=\eI\otimes\eO_{X}$ defines the subscheme $Y_\eJ\subset X$ whose reduction is $(Y,\eO_Y)$. We have the diagram: 
\begin{equation}\label{eq:blup}
\xymatrix@C=4em@R=1.2em{
\tld X:=\Bl_{Y_\eJ}(X)\ar[r]^-{\tld\iota}\ar[d]_-{\si}
\ar@/^4ex/@<+.5ex>[rr]|{\;\phi_X\;}
&\Bl_\eI(\mbb P^N)\ar[r]^-{\phi}\ar[d]^-{B}
&\mbb P^{N'}\ar@{=}[d]
\\ 
X\ar[r]^-\iota&\mbb P^N\ar@{-->}[r]&\mbb P^{N'}
}
\end{equation}

\begin{m-lemma}\label{lm:B}
The diagram \eqref{eq:blup} has the following properties: 
\begin{enumerate}[leftmargin=5ex]
\item[\rm(i)] 
The exceptional divisor of $B$ is $\phi$-relatively ample, hence the exceptional divisor of $\si$ is $\phi_X$-relatively ample. 
\item[\rm(ii)] 
The morphism 
$\phi:\Bl_{\eI}(\mbb P^{N})\to\mbb P^{N'}$ is $G_m$-invariant and 
\begin{equation}\label{eq:estim}
\dim\phi_X(\tld X)=\dim\phi_X(X\sm Y^+)\les\dim (X\sm Y^+).
\end{equation}
\end{enumerate}
\end{m-lemma}

\begin{m-proof}
(i) The subscheme determined by $\eI$ is the vanishing locus of a section in a direct sum of ample line bundles over $\mbb P^N$, so~\ref{prop:p} applies. 

\nit(ii) It holds: 
$\dim\phi_X(\Bl_\eJ(X))=\dim\ovl{\phi_X(X\sm Y)}\;\text{and}\;\;
\phi_X(X\sm Y)=\phi_X(X\sm Y^+)\cup\phi_X(Y^+\sm Y).$ 
For $[\bz_0,\bz']\in Y^+\sm Y$ and $t\in G_m$, the $G_m$-invariance of $\phi_X$ yields:  
$$
\phi_X\big([\bz_0,\bz']\big)=\phi_X\big(t\times[\bz_0,\bz']\big)
=\phi_X\big(\lim_{t\to\infty}t\times[\bz_0,\bz']\big).
$$
But $\disp\lim_{t\to\infty}t\times[\bz_0,\bz']=[0,\bz'']\in X\sm Y^+$, which implies $\phi_X(Y^+\sm Y)\subset\phi_X(X\sm Y^+)$. 
\end{m-proof}

Now we can estimate the ampleness of the source $Y$. 

\begin{m-theorem}\label{thm:y-sink}
Let $X$ be a smooth $G_m$-variety with source $Y$, and $p:=\codim(X\sm Y^+)-1.$ 
The following statements hold: 
\begin{enumerate}[leftmargin=5ex]
\item[\rm(i)] 
The thickening $Y_\eJ$ of $Y$ in \eqref{eq:blup} is a $(\dim Y-p)$-ample subscheme of $X$; in particular, $Y$ is a $p^\apos$ subvariety. 
\item[\rm(ii)]
If $G_m$ acts on $\eN_{Y/X}$ by scalar multiplication, then $Y\subset X$ is $p^\pos$.
\end{enumerate}  
\end{m-theorem}

\begin{m-proof}
(i) We apply the Proposition~\ref{prop:x-b}: $Y_\eJ$ is a $q$-ample subscheme, with
$$
q=1+\dim\phi(\tld X)-\codim_X(Y)\srel{\eqref{eq:estim}}{\les}
1+\dim(X\sm Y^+)-\codim_X(Y).
$$
\nit(ii) In this case we have $Y^+\cong\unbar{\sf N}:= \Spec\big(\Sym^\bullet\eN_{Y/X}^\vee\big)$, cf. \cite[Remark pp. 491]{bb}. Thus $\unbar{\sf N}\subset X$ is open and $G_m$ acts, fibrewise over $Y$, by scalar multiplication. 

The inclusions 
$\unbar{\sf N}\subset\unbar{\sf N}_{\mbb P^N_{\rm source}/\mbb P^N}
=
\{\,[\bz_{N_0},\bz']\mid\bz_{N_0}\neq0\,\}\subset\mbb P^N$ 
are $G_m$-equivariant. The scalar multiplication on the coordinates $\bz'$ exists globally on $\mbb P^N$, so $X\subset\mbb P^N$ is invariant. Hence the exponents $l_\rho$ in \eqref{eq:phi} are all equal one, $\eJ=\eI_Y\subset\eO_X$. 
\end{m-proof}

\begin{m-remark}\label{rmk:cd}
$X\sm Y^+\subset X\sm Y$ is closed, so~\ref{lm:cd-apos} implies $\,\cd(X\sm Y)=\dim(X\sm Y^+).$ 
This simple answer contrasts the elaborate techniques~\cite{ogus,falt-homog}.
\end{m-remark}

\begin{m-example}\label{expl:o-gr} 
\nit(i) 
Endow $W\cong\kk^{w+1}$, $w+1$ even, with a non-degenerate, symmetric bilinear form $\beta$. Let $X:=\oGrs(u+1;W)$ be the orthogonal Grassmannian of $(u+1)$-dimensional isotropic subspaces; in particular, $w+1\ges2(u+1)$. Choose a Lagrangian decomposition $W=\kk^{(w+1)/2}\oplus\kk^{(w+1)/2}$ such that  
$$
\beta=\left[\begin{array}{cc}0&\bone_{(w+1)/2}\\ \bone_{(w+1)/2}&0\end{array}\right],
\quad\text{($\bone$ stands for the identity matrix)}.
$$ 
Consider the action $G_m\srel{\lda}{\to}\SO_{(w+1)/2},\,\lda(t)=\diag\big[t^{-1},\bone_{(w-1)/2},t,\bone_{(w-1)/2}\big],$ whose source is $Y=\{U\mid s:=(1,0,\ldots,0)\in U\}\cong\oGrs(u;w-1)$. For $U\in Y$, $\lda$ acts with weight $t$ on $\eN_{Y/X,U}=\Hom(\lran{s},\lran{s}^\perp/U)$. 
$$
X\sm Y^+=\big\{U\in X\mid s\not\in\underset{t\to 0}{\lim}\lda(t)U\big\}
=
\{U\mid U\subset W':=\kk^{(w-1)/2}\oplus\kk^{(w+1)/2}\}.
$$ 
Let $\lran{s'}:=\Ker(\beta{\rst_{W'}})$: if $w=2u+1$, then $s'\in U$ for all $U\in X\sm Y^+$; for larger $w$, this is not the case. It follows: 
$$
\codim(X\sm Y^+)=\bigg\{ 
\begin{array}{cl}
u&\text{if } w=2u+1;\\ u+1&\text{if }w\ges2u+3,
\end{array}
\;\Rightarrow\;
\text{$Y\subset X$ is: } 
\bigg\{
\begin{array}{cl}
(u-1)^\pos&\text{if }w=2u+1;\\ u^\pos&\text{if }w\ges2u+3.
\end{array}
$$ 

\nit(ii) 
With the previous notation, let $\omega$ be the skew-symmetric bilinear form
$$
\omega=\left[\begin{array}{cc}0&\bone_{(w+1)/2}\\ -\bone_{(w+1)/2}&0\end{array}\right].
$$
Let $X:=\spGrs(u+1;W)$ be the symplectic Grassmannian of $(u+1)$-dimensional isotropic subspaces of $W$. The action of $\,\lda:G_m\to\Sp_{(w+1)/2},$ $\lda(t)=\diag\big[t^{-1},\bone_{(w-1)/2},t,\bone_{(w-1)/2}\big]$ has the source 
$Y=\{U\mid s:=(1,0,\ldots,0)\in U\}\cong\spGrs(u;w-1)$. 

Note that $\eN_{Y/X,U}\cong\Hom(\lran{s},W/U)$, so $G_m$ acts by weight $t^2$ on $\Hom(\lran{s},W/\lran{s}^\perp)$ and weight $t$ on the complement. As before, it holds $\codim(X\sm Y^+)=1+u$. Hence $Y$ is $u^\apos$; more precisely, there is a non-reduced scheme with support $Y$ which is $u^\pos$. 
\end{m-example}


\subsection{Subvarieties of homogeneous varieties}\label{ssct:subvar-homog}

Results due to Faltings, Barth-Larsen, Ogus show that the subvarieties of homogeneous spaces enjoy positivity properties.

\begin{thm-nono}{}
Given a rational homogeneous variety $X=G/P$, with $G$ is semi-simple; let $\ell$ be the minimal rank of its simple factors. Let $Y\subset X$ be a smooth subvariety of codimension $\delta$. The following statements hold: 
\begin{enumerate}[leftmargin=5ex]
\item[\rm(i)]{\rm(cf. \cite[Satz~5, Satz~7]{falt-homog})} $Y$ is $(\ell-2\delta+1)^\pos$. 
\item[\rm(ii)]
$Y\subset\mbb P^n\mbox{ is }p^\pos\;\Leftrightarrow\;
\res_Y^t:H^t(\mbb P^n;\mbb Q)\to H^t(Y;\mbb Q) 
\mbox{ is an isomorphism},\;\forall\,t<p.$
\end{enumerate}
\end{thm-nono}

\begin{m-proof}
(i) $\cd(X\sm Y)\les\dim X-\ell+2\delta-2$ and $\eT_X$ is $(\dim X-\ell)$-ample. Since $\eN_{Y/X}$ is a quotient of $\eT_X$, we conclude by~\ref{prop:N+cd}. 

\nit(ii) $\eN_{Y/\mbb P^n}$ is ample. By \cite[Theorem 4.4, 2.13]{ogus}, $\cd(\mbb P^n\sm Y)<n-p$ if and only if $\res_Y^t$ is an isomorphism and the local cohomological dimension of $Y\subset\mbb P^n$ is at most $n-p$. The latter equals $\codim_{\mbb P^n}Y=n-\dim Y$, since $Y$ is smooth. 
\end{m-proof}
Note that the techniques in~\cite{falt-homog} used for proving (i) above are involved, yet the estimate is not optimal: \ref{cor:Yl} shows that the special Schubert cycles are actually much more positive.


\appendix
\section{Background material}\label{sct:def}

\subsection{Cohomological \textit{q}-ampleness}\label{ssct:cohom-q} 
This notion was introduced by Arapura and Totaro. 

\begin{m-definition}\label{def:q-line}
Let $Y$ be a projective scheme, $\eA\in\Pic(Y)$ an ample line bundle.
\begin{enumerate}[leftmargin=5ex]
\item[\rm(i)]  \cite[Theorem 7.1]{tot} 
An invertible sheaf $\eL$ on $Y$ is \emph{$q$-ample} if, for any coherent sheaf $\eG$ on $X$, holds: 
$\,\exists\,\cst{\eG}\;\forall\,a\ges \cst{\eG}\;\forall\,t>q,\; H^t(Y,\eG\otimes\eL^{a})=0.$
\item[]
It's enough to verify the property for $\eG=\eA^{-k}, k\ges1$ (cf. \cite[Theorem 6.3, 7.1]{tot}).

\item[\rm(ii)]  \cite[Lemma 2.1, 2.3]{arap} 
A locally free sheaf $\eE$ on $Y$ is \emph{$q$-ample} if $\eO_{\mbb P(\eE^\vee)}(1)$ on $\mbb P(\eE^\vee):=\Proj(\Sym^\bullet_{\eO_Y}\eE)$ is $q$-ample. It is equivalent saying that, for any coherent sheaf $\eG$ on $Y$, there is $\cst{\eG}>0$ such that: 
$\,H^t(Y,\eG\otimes\Sym^a(\eE))=0,\;\forall t>q,\; \forall a\ges \cst{\eG}.$ 
\item[] 
The \emph{$q$-amplitude of $\eE$}, denoted by $q^\eE$, is the smallest integer $q$ with this property. 
Note that $\eE$ is $q$-ample if and only if so is $\eE_{Y_{\red}}$ (cf. \cite[Corollary 7.2]{tot}). Also, any locally free quotient $\eF$ of $\eE$ is still $q$-ample; indeed, $\eO_{\mbb P(\eF^\vee)}(1)=\eO_{\mbb P(\eE^\vee)}(1)\otimes\eO_{\mbb P(\eF^\vee)}$.
\end{enumerate}
\end{m-definition}


\subsection{\textit{q}-positivity}

\begin{m-proposition}{\rm \cite[Proposition 1.7]{soms}}\label{prop:q2}
For a \emph{globally generated}, locally free sheaf $\eE$ on $Y$, the following statements are equivalent: 
\begin{enumerate}[leftmargin=5ex]
\item[\rm(i)] 
$\eE$ is $q$-ample (cf. Definition \ref{def:q-line}); 
\item[\rm(ii)] 
The fibres of the morphism $\mbb P(\eE^\vee)\to |\eO_{\mbb P(\eE^\vee)}(1)|$ are at most $q$-dimensional. 
\end{enumerate}
We say that $\eE$ is \emph{Sommese-$q$-ample} if it satisfies any of these conditions.
\end{m-proposition}

\begin{m-definition}\label{def:q-pos}(cf. \cite{andr+grau}) 
Suppose $X$ is a smooth, complex projective variety. A line bundle $\eL$ on $X$ is \emph{$q$-positive}, if it admits a Hermitian metric whose curvature is positive definite on a subspace of $\eT_{X,x}$ of dimension at least $\dim X-q$, for all $x\in X$;  equivalently, the curvature has at each point $x\in X$ at most $q$ negative or zero eigenvalues. 
\end{m-definition}

\begin{m-theorem}\label{thm:mats-bott}\quad
\nit{\rm(i)} \cite[Proposition 28]{andr+grau} $\;q$-positive line bundles are $q$-ample.

\nit{\rm(ii)} {\cite[Theorem 1.4]{mats}} 
Assume $\eE$ is globally generated. Then it holds: 
$$\eE\text{ is Sommese-$q$-ample}\;\Leftrightarrow\;
\eO_{\mbb P(\eE^\vee)}(1)\text{ is $q$-positive.}$$ 

\nit{\rm(iii)} {\cite{bott-lefz,ottm}} 
Let $\eL\in\Pic(X)$ be $q$-positive and $Y\in|\eL|$ a smooth divisor. Then it holds: 
$$
H^t(X;\mbb Z)\to H^t(Y;\mbb Z)\text{ is } 
\bigg\{\begin{array}{rl}
\text{an isomorphism, for}&t\les \dim X-q-2;
\\[1ex] 
\text{injective, for}&t=\dim X-q-1.
\end{array}
$$
\end{m-theorem}


\end{document}